\documentclass[11pt,oneside,reqno]{amsart}
\usepackage{latexsym}
\usepackage{amssymb}
\usepackage{amsthm}
\usepackage{amsmath}
\usepackage{amssymb}
\usepackage[usenames]{color}
\usepackage[mathscr]{euscript}
\usepackage[all]{xy}

\newtheorem{theorem}{Theorem}[section]
\newtheorem{lemma}[theorem]{Lemma}
\newtheorem{corollary}[theorem]{Corollary}
\newtheorem{proposition}[theorem]{Proposition}

\newtheorem{example}[theorem]{Example}
\theoremstyle{definition}
\newtheorem{definition}[theorem]{Definition}

\theoremstyle{remark}
\newtheorem{remark}[theorem]{Remark}

\def\supp{\mathop{\rm supp}\nolimits}
\newcommand{\restr}{\hskip-4pt \restriction}
\newcommand{\Gr}{\mathrm{Gr}}
\newcommand{\gd}{G_{\delta}}

\newcommand{\C}{{\mathbb C}}

\newcommand{\T}{\mathbb T}
\newcommand{\Z}{{\mathbb Z}}

\newcommand{\cont}{\mathfrak c}

\newcommand{\w}{{\wedge}}
\newcommand{\ww}{{\wedge\wedge}}
\newcommand{\la}{\langle\,}
\newcommand{\ra}{\,\rangle}
\newcommand{\rest}{\hskip-2pt\restriction\hskip-2pt}
\newcommand{\restri}{\hskip-3pt\restriction\hskip-3pt}

\def\vart{\vartriangleright}
\newcommand{\p}{{{\text{\tiny$P$}}}}
\newcommand{\h}{{{\text{\tiny$H$}}}}
\def\hull#1{\langle#1\rangle}
\DeclareMathOperator{\ho}{Hom}

\newcommand{\ctg}{\mathcal{T}_{\mbox{\tiny $G$}}}

\author{J.~Galindo} \thanks{Research of the first-listed author supported by
 \emph{Ministerio de Economia y Competitividad}, grant number MTM2011-23118
 and Fundaci\'o Caixa Castell\'o-Bancaixa, grant number P1$\cdot$1B2011-30.}
\address{\noindent Jorge Galindo\\
IMAC, Instituto Universitario de Matem\'aticas y Aplicaciones\\ Universidad
Jaume I, E-12071, Cas\-tell\'on, Spain.\hfill\break \noindent
E-mail: {\tt jgalindo@mat.uji.es}}
\author{M.~Tkachenko}
\address{\noindent M.~Tkachenko, Depto. de Matem\'aticas,
Universidad Aut\'onoma Metropolitana, Av. San Rafael Atlixco 186,
Col. Vicentina, Iztapalapa, C.P. 09340, Mexico D.F., Mexico
\hfill\break \noindent
E-mail: {\tt mich@xanum.uam.mx}}
\author{M.~Bruguera}
\address{\noindent M.~Bruguera, Dept. de Matem\'atica Aplicada I,
Universidad Polit\'ecnica de Catalu\~{n}a, Barcelona, Spain
e-mail: m.montserrat.bruguera@upc.es\hfill\break
\noindent
E-mail: {\tt m.montserrat.bruguera@upc.es}}
\author{C.~Hern\'andez}
\address{\noindent C.~Hern\'andez, Depto. de Matem\'aticas,
Universidad Aut\'onoma Metropolitana, Av. San Rafael Atlixco 186,
Col. Vicentina, Iztapalapa, C.P. 09340, Mexico D.F., Mexico\hfill\break
\noindent E-mail: {\tt chg@xanum.uam.mx}}
\thanks{The second and third listed authors were partially supported
by the CRM (Catalu\~{n}a, Spain) during their visits in May--July, 2008}


\thanks{The major part of this work was done while the first listed author
was visiting the Department of Mathematics of Universidad Aut\'onoma
Metropolitana (UAM), M\'exico in May 2010. Support and hospitality are
gratefully acknowledged.}

\date{March 30, 2012}

\title{Reflexivity in precompact  groups and extensions}

\begin{document}

\begin{abstract}
We establish some general principles and find some
counter-examples concerning the Pontryagin reflexivity
of precompact groups and $P$-groups. We prove
in particular that:
\begin{enumerate}
  \item A precompact Abelian group $G$ of bounded order is reflexive iff
           the dual group $G^\wedge$ has no infinite compact subsets and
           every compact  subset of $G$ is contained in a compact subgroup
           of $G$.
   \item Any extension of a reflexive $P$-group by another reflexive
           $P$-group is again reflexive.
\end{enumerate}
We show on the other hand that an extension of a compact group by
a reflexive $\omega$-bounded group (even dual to a reflexive $P$-group)
can fail to be reflexive.

We also show that the $P$-modification of a reflexive $\sigma$-compact
group can be non-reflexive (even if, as proved in \cite{jlm}, the $P$-modification
of a locally compact Abelian group is always reflexive).
\end{abstract}

\maketitle

\section{Introduction}
The papers \cite{ardaetal,galimaca11,jlm,jlm2} have unveiled that
the duality properties  of the class of precompact groups are more
complicated than expected. The following theorem summarizes some
of the known facts that concern the duality of precompact groups, see
below for unexplained terminology.

\begin{theorem}\label{gendualprec}
Let $G$ be an Abelian group and let $\tau_\h$ be a precompact topology
on $G$ induced by some group of homomorphisms $H\subset \ho(G,\T)$.
The topological group $(G,\tau_\h)^\wedge$ dual to the precompact group
$(G,\tau_\h)$ can be:
  \begin{enumerate}
    \item a discrete group, as for instance when $H$ is countable
     (see \cite{Aub} and \cite{Cha}) or $(G,\tau_\h)$ is the
     $\Sigma$-product of uncountably many copies of the discrete group
     $\Z(2)$ (this can be deduced from (the proof of) Lemma~27.11
     of \cite{bana91}, see Lemma \ref{sigma} below). In this case
     $(G,\tau_\h)$ is not reflexive.
    \item a nondiscrete $P$-group. This is the case  when when $(G,\tau_\h)$
        is the $\omega$-bounded group that arises as the dual of a reflexive
        $P$-group, as those constructed in \cite{jlm} and \cite{jlm2}. Obviously
        $(G,\tau_\h)$ is reflexive in this case.
     \item a precompact noncompact group, as is the case of the infinite
            pseudocompact groups with no infinite compact subsets constructed
            in \cite{galimaca11} and \cite{ardaetal}. These groups are reflexive.
     \item a compact group, as happens when $H=\ho(G,\T)$, the family
        of all homomorphisms of $G$ to $\T$.
  \end{enumerate}
\end{theorem}
The bases on which the reflexivity of precompact groups stands remain
elusive so far. In this paper we give a first insight  to this issue by establishing
some general facts and giving some counterexamples to what could be
regarded as reasonable generalizations of known results. We prove in
Proposition~\ref{charbded} that a precompact Abelian group of bounded
order is reflexive if and only if the compact subsets of the dual group
$G^\wedge$ are finite and every compact subset of $G$ is contained
in a compact subgroup.

We will especially address the behavior of reflexivity under
extensions in precompact, $\omega$-bounded, and $P$-groups.
We recall that the class of $P$-groups is naturally linked to that of
precompact groups through duality since, by \cite[Lemma~4.1]{jlm},
the dual group $G^\wedge$ of every $P$-group $G$ is $\omega$-bounded
and hence precompact.

Our starting point in this regard is the fact that an extension of a reflexive
group by a compact group is again reflexive (see \cite[Theorem~2.6]{BCM}).
We show in Example~\ref{ex:extref} that an extension of a compact group
by a reflexive precompact (even $\omega$-bounded) group may be
non-reflexive. In Corollary~\ref{cor:ext2} we prove, in contrast, that an
extension of a reflexive $P$-group by another reflexive $P$-group is
always reflexive.

Answering a question in \cite{jlm} we prove in Section~\ref{Sec:Pmod}
that, while the $P$-modification of an LCA group is always reflexive, the
$P$-modification of a reflexive $\sigma$-compact group can fail to be
reflexive.

\subsection{Notation and terminology}
All groups considered here are assumed to be Abelian.
A \textit{character\/} of a topological group $G$ is a continuous
homomorphism of $G$ to the circle group $\T=\{z\in\C: |z|=1\}$ when
the latter is considered as a subgroup of in the complex plane $\C$
with its usual topology and multiplication. The group $G^\wedge$ of
all characters of $G$ with the pointwise multiplication is called the
\textit{dual group} or simply the \textit{dual} of $G$. The dual group
$G^\wedge$ carries the \textit{compact-open} topology $\tau_{co}$
defined as follows.

Put $\T_+=\{z\in\T: Re(z)\geq 0\}$.  For a nonempty set $K\subset G$,
we define
\[
K^\vartriangleright=\left\{\chi\in G^\wedge\colon \chi(K)\subset
\T_+    \right\}.
\]
The collection of sets $\left\{K^\vartriangleright\colon K\subset
G,\ K \mbox{ is compact} \right\}$ forms a neighborhood basis at the
identity of $G^\wedge$ for the compact-open topology $\tau_{co}$.
Let us note that the sets $K^\vartriangleright$ are not necessarily
open in $(G^\wedge,\tau_{co})$ since $\T_+$ is not open in $\T$.
If, instead of $\T_+$, we use a smaller neighborhood $U$ of $1$
in $\T$ to construct the sets $K^\vartriangleright$, the resulting sets will
also form a neighbourhood basis at the identity  for the compact-open
topology $\tau_{co}$.

The subgroup
\[
K^\perp=\left\{\chi\in G^\wedge\colon\chi(K) =\{1\}\right\}
\]
of $G^\wedge$ is called the \textit{annihilator} of a set $K\subset G$.
If $B\subset G^\wedge$, we will also find useful to refer to the set
$B^\vartriangleleft=\left\{x\in G\colon \chi(x)\in \T_+ \mbox{ for all }
\chi \in B\right\}$.

A subset $A$ of a topological group $G$ is called \textit{quasi-convex}
if $A=(A^\vartriangleright)^\vartriangleleft$. If quasi-convex sets in $G$
form a neighborhood base of the neutral element in $G$, the group $G$
will be called \textit{locally quasi-convex} \cite{Aub}.

The \textit{bidual group} of $G$ is $G^{\wedge\wedge}=(G^\wedge)^\wedge$.
The \textit{evaluation mapping} $\alpha_G\colon G\to G^{\wedge\wedge}$
is defined by
\[
\alpha_G(x)(\chi)=\chi(x),
\]
for all $x\in G$ and $\chi\in G^\wedge$. It is easy to see that $\alpha_G$
is a homomorphism. If it is a topological isomorphism of $G$ onto
$G^{\wedge\wedge}$, the group $G$ is called \textit{reflexive.} Every
reflexive group is locally quasi-convex \cite[Prop.~6.6]{Aub}.

Let $H$ be a subgroup of a topological group $G$. We say that $H$
{\it dually embedded\/} in $G$ if each continuous character of $H$ can
be extended to a continuous character of $G$. A subgroup $H$ of a
topological group $G$ is said to be {\it$h$-embedded\/} into $G$
provided that any homomorphism $\varphi$ of $H$ to an arbitrary
compact group $K$ is extendable to a continuous homomorphism
$\tilde\varphi\colon G\to K$. Note that if $H$ is an $h$-embedded
subgroup of $G$, then any homomorphism of $H$ to a compact group
is continuous. Note also that every $h$-embedded subgroup $H$ of
$G$ is dually embedded in $G$.

Let $N$ be a closed subgroup of a topological group $G$.
The group $G$ is usually called an \textit{extension} of $G/N$ by
the group $N$. For example, every feathered (equivalently, almost
metrizable) Abelian group is an extension of a metrizable group by
a compact group (see \cite[Theorem~4.3.20]{AT}).

A space $X$ is called \textit{$\omega$-bounded} provided the closure
of every countable subset of $X$ is compact. It is clear that every
$\omega$-bounded space is countably compact, but not vice versa.

Every pseudocompact topological group is precompact
\cite[Theorem~1.1]{CR2}. Hence all $\omega$-bounded and countably
compact groups are precompact. By a well-known theorem of Comfort
and Ross in \cite{CR}, a topological group $(G,\tau)$ is precompact
if and only if the topology $\tau$ of $G$ is the topology $\tau_\h$
generated by a group of characters $H\subset \ho(G,\T)$.

By a \textit{protodiscrete} group we understand a topological group
having a basis of neighborhoods of the identity consisting of open
subgroups. Protodiscrete groups are also known as linear groups.
Evidently, protodiscrete Abelian groups are locally quasi-convex.

A \textit{$P$-space} is a space in which every $G_\delta$-set is open.
A \textit{$P$-group} is a topological group which is a $P$-space.
According to \cite[Lemma~4.4.1]{AT}, every $P$-group is protodiscrete.

Given a space $X$, the \textit{$P$-modification} of $X$, denoted by
$PX$, is the underlying set $X$ endowed with the topology whose base
consists of $G_\delta$-sets in the original space $X$. It is clear that
the $P$-modification of a topological group is again a topological group.

Let $A$ denote a non-empty index set and let, for each $\alpha \in A$,
$G_\alpha$ be a compact group with identity  $e_\alpha$.
Given  $x\in \prod_{\alpha \in A}G_\alpha$, we put
\[
\supp(x)=\{\alpha\in A: x_\alpha \neq e_\alpha\},
\]
\[
\sum\prod_{\alpha \in A}G_\alpha=\{x\in\prod_{\alpha \in A}G_\alpha
\colon |\supp(x)|\leq\omega\},
\]
and
\[\bigoplus_{\alpha \in A}G_\alpha=\{x\in\prod_{\alpha \in A}G_\alpha
\colon |\supp(x)|<\omega\}.
\]
It is clear that $\sum\prod_{\alpha \in A}G_\alpha$ is a dense
$\omega$-bounded subgroup of $\prod_{\alpha \in A}G_\alpha$
(see \cite[Proposition~1.6.30]{AT}). In particular, the group
$\sum\prod_{\alpha \in A}G_\alpha$ is countably compact.
This subgroup of $\prod_{\alpha \in A}G_\alpha$ is called the
 \textit{$\Sigma$-product} of the family $\{G_\alpha\colon
 \alpha \in A\}$. When $G_\alpha=G$ for all $\alpha \in A$, we use
 the symbol $\Sigma G^A$ instead of $\sum\prod_{\alpha \in A}G_\alpha$.
The group $\oplus_{\alpha\in A} G_\alpha$ is usually known as
the \textit{direct sum} (and also as the \textit{$\sigma$-product})
of the family $\{G_\alpha\colon \alpha \in A\}$.

\section{On the duality of precompact Abelian groups}
In this section we collect several general results concerning the
duality theory of precompact groups. Some of them appear in the
literature, while others might be known to specialists, but it seems
convenient to have them all collected here for future references.
In any case, the two lemmas below are well known.

\begin{lemma}\label{Le:sepa}
If the family $G^\wedge$ separates elements of the group $G$, then
$\alpha_G\colon G\to G^{\wedge\wedge}$ is a monomorphism.
\end{lemma}

\begin{proof}
If $g\in G$ and $g\neq 0_G$, then there exists $\chi\in G^\wedge$
such that $\chi(g)\neq 1$ or, equivalently, $\alpha_G(g)(\chi)\neq 1$.
Hence $\alpha_G(g)$ is distinct from the neutral element of
$G^{\wedge\wedge}$ and $\alpha_G$ is a monomorphism.
\end{proof}

\begin{lemma}\label{duaemb}
Every subgroup $L$ of a precompact group $G$ is dually embedded.
\end{lemma}

\begin{proof}
Let $\varrho G$ be the Weil completion of $G$ and $H$ be the closure
of $L$ in $\varrho G$. Then $H$ is a closed subgroup of the compact
Abelian group $\varrho G$, so \cite[Prop.~9.6.2]{AT} implies that $H$
is dually embedded in $\varrho G$. By \cite[Prop.~3.6.12]{AT}, the dense
subgroup $L$ of $H$ is dually embedded in $H$. Hence $L$ is dually
embedded in $\varrho G$ and in $G$.
\end{proof}

\begin{proposition}\label{compinref}
Let $G$ be a topological group. If each compact subset of $G$ is contained
in a reflexive, dually embedded subgroup, then $\alpha_G$ is onto.
\end{proposition}

\begin{proof}
Let $\Psi\in G^\ww$ be an element of the bidual group. Then there is a compact
subset $K\subset G$ such that $K^\vartriangleright\subset \Psi^{-1}(\T_+)$.
By hypothesis there is a reflexive dually embedded subgroup $L$ of $G$ with
$K\subset L$. Since $L^\perp=L^\vartriangleright\subset K^\vartriangleright$
we have that $L^\perp \subset \ker \Psi$. Since $L$ is dually embedded,   this implies that $\Psi$ factorizes
through $\widehat{L}$. In other words, there is a continuous homomorphism
$\Psi_L\colon \widehat{L}\to\T$ that, denoting by $R$ the canonical restriction
map of $G^\wedge$ to $L^\wedge$, makes the following
diagram commute.
\[
\xymatrix{ G^\wedge \ar@{>}[rr]^R \ar[dr]^{\Psi} & & L^\wedge
\ar[dl]_{\Psi_L} \\ & \T &}
\]
Notice that $R$ is a continuous surjective homomorphism, while the
continuity of $\Psi_L$ follows from the inclusion $K\subset L$. Since
$L$ is reflexive, $\Psi_L=\alpha_L(g)$ for some $g\in L$. This means
that $\Psi=\alpha_G(g)$. Indeed, if $\chi\in G^\wedge$ then $\Psi(\chi)=
\Psi_L(\chi\restr_L)= \chi\restr_L(g)=\chi(g)$.
\end{proof}

The applicability of Proposition~\ref{compinref} in our context is enhanced by
the following theorem of S.~Hern\'andez and S.~Macario.

\begin{lemma}[\cite{hernmaca03}]\label{pseudo}
If $G$ is a pseudocompact group, then $G^\wedge$ has no infinite
compact subsets.
\end{lemma}

\begin{proof}
Let $H$ be a subgroup of $\ho(G,\T)$ that induces the topology
of $G$. By \cite[Proposition~3.4]{hernmaca03}, every countable
subgroup of $(H,\ctg)$ is $h$-embedded, where $\ctg$ is the topology
(of $H$) of pointwise convergence on elements of $G$. This implies
that $(H,\ctg)$ cannot contain infinite compact subsets, see
\cite[Prop.~2.1]{ardaetal} or \cite{galimaca11}. The same is true,
\emph{a fortiori}, for $G^\wedge$.
\end{proof}

\begin{corollary}
\label{caromega} Let $G$ be a pseudocompact group. If every compact
subset of $G$ is contained in a compact subgroup of $G$, then $G$
is reflexive.
\end{corollary}

\begin{proof}
Since $G$ is pseudocompact (hence precompact), $G^\wedge$
separates elements of $G$. Therefore, Lemma~\ref{Le:sepa} shows
that $\alpha_G$ is a monomorphism. By Lemma~\ref{duaemb},
all subgroups of  $G$ are dually embedded. Hence Proposition~\ref{compinref}
implies that $\alpha_G$ is a group isomorphism.  Since, by
Lemma~\ref{pseudo}, $G^\ww$ carries the topology of pointwise
convergence on elements of $G^\wedge$, just as $G$ does, this
isomorphism is easily seen to be a homeomorphism.
The group $G$ is therefore reflexive.
\end{proof}

\begin{proposition}\label{subref}
If $G$ is precompact and $\alpha_G$ is onto, then every closed metrizable
subgroup of $G$ is compact.
\end{proposition}

\begin{proof}
Let $N$ be a metrizable subgroup of $G$. If $\varrho{N}$ denotes
the completion of $N$ we have, as a consequence of the
Au{\ss}enhofer--Chasco theorem (see \cite{Aub} or \cite{Cha}) that
$N^\ww =\varrho{N}$.

If some closed metrizable subgroup $N$ of $G$ is not compact,
there is $\Psi\in N^\ww$ with $\Psi\notin \alpha_N(N)$. Define now
$\Psi^G\colon G^\wedge \to \T$ by $\Psi^G(\chi)=\Psi(\chi\restr_N)$
for each $\chi\in G^\wedge$. Then $\Psi^G\in G^\ww$. Since $\alpha_G$
is onto by hypothesis, there must be  $g\in G$ with $\Psi^G=\alpha_G(g)$.
As $N$ is closed it follows that $g\in N$, for otherwise there is
$\chi\in G^\wedge$ with $\chi\restr_N\,=1$ and $\chi(g)\neq 1$ yielding
that $\chi(g)=\alpha_G(g)(\chi)= \Psi^G(\chi)=\Psi(\chi\restr_N)=1$. But, since
$g\in N$, this implies $\Psi=\alpha_N(g)\in\alpha_N(N)$, against our choice
of $\Psi$.
\end{proof}

\begin{proposition}\label{prop:charprec}
Let $G$ be a precompact topological group such that $G^\wedge$ is
a protodiscrete group with no infinite compact subsets. Then $G$ is reflexive
if and only if every compact subset of $G$ is contained in  a compact
subgroup of $G$.
\end{proposition}

\begin{proof}
Observe that the general hypothesis implies that $\alpha_G$ is a
topological isomorphism of $G$ onto a subgroup of $G^{\wedge\wedge}$.

  \emph{Sufficiency:} If every compact subset of $G$ is
contained in  a compact subgroup of $G$, then $\alpha_G$ is onto
by Lemma~\ref{compinref}.

  \emph{Necessity:} Suppose that $G$ is reflexive and let $K$ be a compact
subset of $G$. Since  $K^\vartriangleright$ is a neighbourhood of the
identity of $G^\w$, there is an open subgroup $N$ of $G^\wedge$ with
$N\subset  K^\vartriangleright$. Now
$\alpha_G(K)\subset K^{\vartriangleright\vartriangleright}\subset N^\perp$
and $N^\perp$ is a compact subgroup of $G^\ww=G$, for $N$ is open in
$G^\wedge$ (see \cite[Lemma~2.2]{Nob2}). Thus $K$ is contained in the
compact subgroup $\alpha_G^{-1}(N^\perp)$ of $G$.
\end{proof}

\begin{corollary}
Let $G$ be a precompact topological group such that $G^\wedge$ is
a $P$-group. Then $G$ is reflexive if and only if every compact subset
of $G$ is contained in a compact subgroup of $G$.
\end{corollary}

\begin{proof}
Clearly all compact subsets of a $P$-group are finite. Since every
$P$-group is protodiscrete \cite[Lemma~4.4.1]{AT}, the conclusion
follows from Proposition~\ref{prop:charprec}.
\end{proof}

\begin{lemma}\label{lem:bded}
If $G$ is a precompact group of bounded order, then $G^\wedge$ is a
protodiscrete group.
\end{lemma}

\begin{proof}
Let $n$ be the exponent of $G$. Choose a neighbourhood $V$ of 1
in $\T$ not containing $n$-roots of 1 other than 1 itself. Then the
equality $\left\{\chi\in G^\wedge\colon \chi(K)\subset V\right\}=K^\perp$
holds for every compact set $K\subset G$. These sets form a basis of
open neighbourhoods of the identity in $G^\wedge$ and each of
them is evidently a subgroup of $G^\wedge$.
\end{proof}

\begin{proposition}\label{charbded}
Let $G$ be a precompact group of bounded order. Then $G$ is
reflexive if and only if it has the following two properties:
\begin{enumerate}
    \item $G^\wedge$ has no infinite compact subsets.
    \item Every compact subset of $G$ is
  contained in  a compact subgroup of $G$.
  \end{enumerate}
\end{proposition}

\begin{proof}
If $G$ is reflexive and $K\subset G^\wedge$ is compact, then $K$ is
finite by \cite[Proposition~2.7]{ardaetal}. By Lemma~\ref{lem:bded},
$G^\wedge$ is protodiscrete. Hence (2) follows from
Proposition~\ref{prop:charprec}.

Suppose conversely that (1) and (2) hold. Then, by Lemma~\ref{lem:bded},
$G$ satisfies the hypothesis of Proposition~\ref{prop:charprec} and is
therefore reflexive.
\end{proof}

Combining Lemma~\ref{pseudo} and Proposition~\ref{charbded} we
deduce the following:

\begin{corollary}\label{cor:boundedpseudo}
A pseudocompact group $G$ of bounded order is reflexive if and only if
each compact subset of $G$ is contained in a compact subgroup of $G$.
\end{corollary}

\begin{remark}
Proposition~\ref{charbded} is false for groups of infinite exponent.
It suffices to consider a torsion-free pseudocompact group $G$ without
infinite compact subsets (see Theorems~5.5,~5.7 and Corollary~5.6
of \cite{galimaca11}). Such a group cannot contain any nontrivial
compact  subgroup.
\end{remark}

\section{Some extension results}\label{Sec:1}
We say that $\mathcal P$ is a \textit{three space property} if for
every Hausdorff topological group $G$ and every closed invariant
subgroup $N$ of $G$ such that both $N$ and $G/N$ have $\mathcal P$,
the group $G$ has $\mathcal P$ as well. It is known, on one hand,
that compactness, connectedness, precompactness,
pseudocompactness, Ra\u{\i}kov completeness, etc., are all three
space properties. On the other hand, Lindel\"ofness, normality, having
a countable network, countable compactness, and many others, fail
to be three space properties (see \cite{BT2} for a detailed discussion
on the subject).

We have already mentioned in the introduction that extensions of
reflexive groups by compact groups preserve reflexivity. However, we
will see in Section~\ref{sec:nonexten} that an extension of a compact
group by a reflexive (even $\omega$-bounded) group may fail to be
reflexive. Therefore reflexivity is not a three space property even among
$\omega$-bounded groups. The situation is completely different for
the class of $P$-groups as we now set on to show. First we need
two auxiliary facts.

\begin{lemma}\label{le:extprot}
Let $H$ be a closed subgroup of a topological Abelian group $G$.
If the groups $H$ and $G/H$ are protodiscrete, so is $G$.
\end{lemma}

\begin{proof}
Take an arbitrary neighborhood $U$ of the neutral element $e$ in
$G$ and choose another neighborhood $V$ of $e$ with $V+V\subset U$.
Since $H$ is protodiscrete, there exists an open subgroup $C$
of $H$ with $C\subset H\cap V$. Let $W$ be an open symmetric
neighborhood of $e$ in $G$ such that $W\subset V$ and
$(W+W+W)\cap H\subset C$.

Denote by $\pi$ the quotient homomorphism of $G$ onto $G/H$.
Since $G/H$ is protodiscrete, one can find an open subgroup $K$
of $G/H$ satisfying $K\subset \pi(W)$. We claim that $N=\pi^{-1}(K)
\cap (W+C)$ is an open subgroup of $G$ with $N\subset U$. The set
$N$ is open in $G$ since $K$ is open in $G/H$ and $W$ is open in $G$.
It is also clear that $N\subset W+C\subset V+V\subset U$. Therefore,
to finish the proof, it suffices to verify that $N$ is a subgroup of $G$.

First we note that $N$ is symmetric, and $C\subset N$ (observe that
$C\subset H$). In fact, our definition of $N$ implies that $N+C=N$.
We then take arbitrary elements $x_1,x_2\in N$. There exist
$w_1,w_2\in W$ and $c_1,c_2\in C$ such that $x_i=w_i+c_i$ for $i=1,2$.
Note that $\pi(N)=K\cap \pi(W+C)=K\cap \pi(W)=K$. Since $K$ is a
subgroup of $G/H$, we see that
$\pi(x_1+x_2)=\pi(x_1)+\pi(x_2)\in K$. Hence there exists $x_3\in N$
such that $\pi(x_3)=\pi(x_1+x_2)$ and we can find $h\in H$ such that
$x_1+x_2=x_3+h$. Choose $w_3\in W$ and $c_3\in C$ such that
$x_3=w_3+c_3$. Then $w_1+w_2-w_3=h-c_1-c_2+c_3$, so the element
$h-c_1-c_2+c_3$ is in $(W+W+W)\cap H\subset C$. In its turn, this implies
that $h\in C$. Therefore, $x_1+x_2=x_3+h\in N+C=N$. This shows that
$N$ is a subgroup of $G$ and finishes the proof.
\end{proof}

\begin{proposition}\label{prop:ext1}
Let $G$ be a topological group and $H$ be a closed subgroup of
$G$. Suppose that both $H$ and $G/H$ are reflexive, protodiscrete, and
contain no infinite compact subsets. If $K\subset G^\wedge$ is compact,
then there is an open subgroup $N$ of $G$ with $K\subset N^\perp$.
\end{proposition}

\begin{proof}
Let $K$ be a compact subset of $G^\wedge$. Denote by $j$ the natural
embedding of $H$ to $G$. We consider the exact sequences
\[
\xymatrix{ H\ar@{>}[rr]^j  &&G\ar@{>}[rr]^\pi && G/H}\]
and
\[
\xymatrix{ \left(G/H\right)^\wedge\ar@{>}[rr]^{\widehat{\pi}}  && G^\wedge
\ar@{>}[rr]^{\widehat{\jmath}} && H^\wedge.}
\]
Then $\widehat{\jmath}(K)$ is a compact subset of $H^\wedge$. Since $H$
is reflexive and protodiscrete, there is an open subgroup $P$ of $H$  such that
$\widehat{\jmath}(K)\subset P^\perp_H$, where $P^\perp_H$ is the annihilator
of $P$ in $H^\wedge$. This implies that $K\subset P^\perp_G$, with $P^\perp_G$
the annihilator of $P$ in $G^\wedge$.

Let now $\widetilde{P}$ be an open subgroup of $G$ with
$\widetilde{P}\cap H=P$ (here we use the protodiscreteness of $G$ provided by  Lemma  \ref{le:extprot})
and consider the sequence of mappings
\[
\xymatrix{
\widetilde{P} \ar@{>}[rr]^{\pi_\p}& &(\widetilde{P}+H)/H
\ar@{>}[rr]^{\sigma}&& \widetilde{P}/P.}
\]
Here $\pi_\p$ is the quotient homomorphism and $\sigma$ is the obvious
group isomorphism given by the third isomorphism theorem. Since
$\widetilde{P}$ is open in $G$, $\sigma$ is a topological isomorphism.

We have the following dual sequence:
\[
\xymatrix{
\left( \widetilde{P}/P\right)^\wedge
\ar@{>}[rr]^{\widehat{\sigma}} & & \left((\widetilde{P}+H)/H\right)^\wedge
 \ar@{>}[rr]^{\widehat{\pi}_{_\p}}& & \widetilde{P}^\wedge}.
 \]
As neither $\widetilde{P}$ nor $(\widetilde{P}+H)/H$ have infinite
compact subsets, $\widehat{\pi}_\p$ (and hence also
$\widehat{\pi}_\p\circ\widehat\sigma$) is a topological isomorphism
onto $P^\perp$, the annihilator of $P$ in $\widetilde{P}^\wedge$.

Finally, denote by $R\colon G^\wedge\to \widetilde{P}^\wedge$
the restriction homomorphism dual to the inclusion of $\widetilde{P}$
into $G$. Since $K\subset P^\perp_G$, $R(K)$ is contained in the
image of $\widehat{\pi}_\p$. Then $\widehat{\pi}_\p^{-1}(R(K))$ is
a compact subset of the dual of the group $(\widetilde{P}+H)/H$.
Notice that $(\widetilde{P}+H)/H$ is reflexive as an open subgroup
of the reflexive group $G/H$ \cite[Theorem~2.3]{BCM}. Since the
former group is protodiscrete, there is an open subgroup $N$ of
$\widetilde{P}$ (hence of $G$) such that $\widehat{\pi}_\p^{-1}(R(K))
\subset \pi_P(N)^\perp$.

We finally claim that $K\subset N^\perp$. To that end, let $\chi\in K$ and
$x\in N$ be arbitrary elements. Since $\widehat{\jmath}(K)\subset P^\perp_H$,
there is some $\widetilde{\chi}\in \left( \widetilde{P}/P\right)^\wedge$ with
$R(\chi)=(\sigma\circ \pi_{\p})^\wedge(\widetilde{\chi})$.
But then $\widehat{\sigma}(\widetilde{\chi})\in\widehat{\pi}_\p^{-1}(R(K))$
and, recalling that $\widehat{\pi}_\p^{-1}(R(K))\subset \pi_P(N)^\perp$,
we see that
\begin{equation}\label{seq}
\widehat{\sigma}(\widetilde{\chi})(\pi_\p(x))=1.
\end{equation}

In addition,

\begin{equation}\label{seq2}
\widehat{\sigma}(\widetilde{\chi})(\pi_\p(x))=
\widehat{\pi}_\p\bigl(\widehat{\sigma}(\widetilde{\chi})\bigr)(x)=
(\sigma\circ \pi_\p)^\wedge(\widetilde{\chi})(x)=R(\chi)(x)=\chi(x).
\end{equation}
Equalities~\eqref{seq} and~\eqref{seq2} show that $\chi(x)=1$
for all $x\in N$ and $\chi \in K$. Hence $K\subset N^\perp$.
\end{proof}

\begin{corollary}\label{cor:ext1}
Let $H$ be a closed subgroup of a topological group $G$ and assume
that $H$ and $G/H$ contain no infinite compact subsets. If both $H$
and $G/H$ are reflexive and protodiscrete, then $G$ is reflexive and
protodiscrete.
\end{corollary}

\begin{proof}
It follows from our assumptions about $H$ and $G/H$ that all compact subsets
of the group $G$ are finite, while Lemma~\ref{le:extprot} implies that $G$ is
protodiscrete and, therefore, locally quasi-convex.

Take a basic open neighborhood of the neutral element in $G^{\wedge\wedge}$
of the form $K^\vartriangleright$, where $K$ is a compact subset of $G^\wedge$.
By Proposition~\ref{prop:ext1}, there exists an open subgroup $N$ of $G$
such that $K\subset N^\perp$. Then $\alpha_G(N)\subset
(N^\perp)^\vartriangleright\subset K^\vartriangleright$, so $\alpha_G$
is continuous.

Since $G$ is locally quasi-convex, $\alpha_G$ is necessarily
injective and open as a mapping onto $\alpha_G(G)$
\cite[Prop.~6.10]{Aub}. Finally, $G^\wedge$ carries the topology of
pointwise convergence on elements of $G$ since $G$ does not have
infinite compact subsets and, as a consequence, $\alpha_G$ is
surjective (apply \cite[Theorem~1.3]{CR}). We conclude that $G$ is
reflexive.
\end{proof}

In Corollary~\ref{cor:ext2} below we show that the property of being
a reflexive $P$-group is closed under extensions. Let us first
establish that the class of $P$-groups behaves similarly:

\begin{lemma}\label{Le:Pg}
Suppose that $H$ is a closed subgroup of a topological (not
necessarily Abelian) group $G$ and that both $H$ and $G/H$ are
$P$-groups. Then $G$ is also a $P$-group.
\end{lemma}

\begin{proof}
Let $\pi\colon G\to G/H$ be the quotient homomorphism. Denote by
$\tau$ the topology of $G$ and let $\tilde\tau$ be the $P$-modification
of the topology $\tau$. It is clear that $\tilde\tau$ is finer than $\tau$
and that $\widetilde{G}=(G,\tilde\tau)$ is again a topological group.
In particular, $H$ is closed in $\widetilde{G}$. Since $H$ is a $P$-group,
we see that $\tilde{\tau}\rest H=\tau\rest H$. Similarly, since $G/H$ is
a $P$-group, the quotient groups $\widetilde{G}/H$ and $G/H$ carry
the same topology, i.e., $\tilde{\tau}/H=\tau/H$. By \cite[Lemma~1]{DS},
this implies that $\tilde{\tau}=\tau$, so $G$ is a $P$-group.
\end{proof}

\begin{corollary}\label{cor:ext2}
If $H$ is a closed subgroup of a topological group $G$ and both $H$
and $G/H$ are reflexive $P$-groups, then so is $G$.
\end{corollary}

\begin{proof}
It follows from Lemma~\ref{Le:Pg} that $G$ is a $P$-group. Notice that
$P$-groups are protodiscrete \cite[4.4.1\,(a)]{AT} and have no infinite
compact subsets \cite[4.K.2]{gilljeri}. Hence $G$ is reflexive by
Corollary~\ref{cor:ext1}.
\end{proof}

\section{Extending compact groups by reflexive groups}\label{sec:nonexten}
In this section we present several examples of non-reflexive
extensions of compact groups by $\omega$-bounded groups.
Let us start with two lemmas.

\begin{lemma}\label{lem:constG}
Let $\phi \colon H_1\to H_2$ be a continuous homomorphism
of topological groups and $Gr(\phi)=\{(x,\phi(x)): x\in H_1\}$ be
the graph of $\phi$ considered as a subgroup of $H_1\times H_2$.
If $L$ is a dense subgroup of $H_2$, then the subgroup $G=\Gr(\phi)+ N$
of $H_1\times H_2$ is an extension of $H_1$ by $L$, where
$N=\{e\}\times L$ and $e$ is the neutral element of $H_1$. If $L$ is
$\gd$-dense in $H_2$, then $G$ is $\gd$-dense in $H_1\times H_2$.
\end{lemma}

\begin{proof}
It is easy to see that $N=G\cap (\{e\}\times H_2)$, so $N$ is a closed
subgroup of $G$ topologically isomorphic to $L$. Since $N$ is dense
in $\{e\}\times H_2$, it follows from \cite[Lemma~1.3]{Gra} that the
restriction to $G$ of the projection $\pi\colon H_1\times H_2\to H_1$ to
the first factor is an open homomorphism of $G$ onto $H_1$.
Therefore, $G/N\cong H_1$.

Finally, suppose that $L$ is a $\gd$-dense subgroup of $H_2$ and
let $V_1\times V_2$ be a non-empty $\gd$-subset of $H_1\times H_2$.
Take $x\in V_1$; by the $\gd$-density of $L$ there must be
$y\in L\cap(\phi(x)^{-1}+V_2)$. Then $(x,\phi(x)+y)\in G\cap (V_1\times V_2)$
and hence $G$ is $\gd$-dense in $H_1\times H_2$.
\end{proof}

\begin{lemma}\label{pdualsigma}
Let $P=P\Z(2)^{\cont}$ be the $P$-modification of the compact
group $\Z(2)^\cont$. There is a topological monomorphism
$j\colon P^\wedge \to \Z(2)^{2^\cont}$ such that
$j(P^\wedge)\cap \Sigma \Z(2)^{2^\cont}=\{e\}$.
\end{lemma}

\begin{proof}
Let $A$ be a maximal independent subset of $P$. Denote by
$j_A$ the restriction homomorphism of $P^\wedge$ to $\Z(2)^A$,
$j_A(\chi)=\chi\restri A$ for each $\chi\in P^\wedge$. We claim that
$j_A$ is a topological isomorphism of $P^\wedge$ onto the dense
subgroup $j_A(P^\wedge)$ of $\Z(2)^A$.

Indeed, $j_A$ is a monomorphism since $A$ generates the group $P$
algebraically. Since the compact subsets of $P$ are finite, $P^\wedge$
is a topological subgroup of $\Z(2)^P$. Hence the continuity of $j_A$
follows from the continuity of the projection of $\Z(2)^P$ to $\Z(2)^A$.
Let us show that $j_A$ is open as a mapping of $P^\wedge$ onto
the subgroup $j_A(P^\wedge)$ of $\Z(2)^A$. Given a neighborhood
$U$ of the neutral element in $P^\wedge$, we can find elements
$x_1,\ldots,x_n$ in $P$ such that
$$
V=\{\chi\in P^\wedge: \chi(x_i)=1\mbox{ for each } i=1,\ldots,n\}\subset U.
$$
For every $i\leq n$, take elements $a_{i,1},\ldots,a_{i,k_i}\in A$
such that $x_i=a_{i,1}\cdot \cdots \cdot a_{i,k_i}$ and let
$$
B=\{a_{i,j}: 1\leq i\leq n,\ 1\leq j\leq k_i\}
$$
and
$$
O=\{x\in\Z(2)^A: x(b)=1 \mbox{ for each } b\in B\}.
$$
An easy verification shows that $j_A(V)\supset O\cap j_A(P^\wedge)$,
which implies that $j_A$ is open. Summing up, $j_A$ is a topological
monomorphism. The density of $j_A(P^\wedge)$ in $\Z(2)^A$ is
evident. Note that $j_A(\chi)\in \Sigma\Z(2)^A$ if and only if
$\chi(a)=1$ for all but countably many $a\in A$.

For every $C\in [\cont]^{\leq \omega}$, let $N_C=\pi_{C}^{-1}(e_C)$
be the $\gd$-subset of $\Z(2)^{\cont}$, where
$\pi_C\colon \Z(2)^{\cont}\to\Z(2)^C$ is the projection and $e_C$
is the neutral element of $\Z(2)^C$. Let $Y_C$
denote a maximal independent subset of $N_C$. Observe that,
for every $x\notin N_C$,  $x\cdot Y_C$ is again independent.

For every $B\in [\cont]^{\leq \omega}$, choose an element
$x_B\in \Z(2)^{\cont}$ which is supported precisely on $B$, that is,
$x_B(\alpha)=-1$ iff $\alpha\in B$.

We then define, for all $B,C\in [\cont]^{\leq \omega}$ with
$B\cap C\neq \emptyset$, the sets  $X_{B,C}=x_B\cdot Y_C$. Since
each $X_{B,C}$ is independent and has cardinality $2^\cont$, we may
construct  as in  Lemma~4.4 of \cite{comfgali03}  a collection of sets
$Z_{B,C}\subset X_{B,C}$, such that:
\begin{enumerate}
\item[(1)]  $|Z_{B,C}|=2^\cont$;
\item[(2)]  $Z_{B,C}\cap Z_{B^\prime,C^\prime}=\emptyset$ if
                $(B,C)\neq ( B^\prime,C^\prime)$;
\item[(3)] $\bigcup_{B,C} Z_{B,C}$ is independent.
\end{enumerate}
Let $Z$ be a maximal independent subset of $\Z(2)^{\cont}$
containing the union in (3). Denote by $j=j_Z$ the topological
monomorphism of $P^\wedge$ to $\Z(2)^Z$ corresponding to $Z$.

Suppose now that $\psi\in P^\wedge$, $\psi\neq\mathbf{1}$. Since
$\psi$ is continuous there must be $C\in [\cont]^{\leq \omega}$
such that $\psi(N_C)=\{1\}$. As $\psi\neq \mathbf{1}$ (and noting
that $\Sigma \Z(2)^\cont$ is dense in $P$), there must also be
$B\in [\cont]^{\leq \omega}$ with $\psi(x_B)=-1$; observe that
necessarily $B\cap C\neq \emptyset$. It follows that
$\psi(x_B\cdot  N_C)=\{-1\}$ and $\psi(Z_{B,C})=\{-1\}$, which means
that $j(\psi)\notin \Sigma\Z(2)^Z$. This implies that
$j(P^\wedge)\cap \Sigma\Z(2)^Z=\{e\}$. Note that $|Z|=2^\cont$,
so we complete the proof by identifying $Z$ with $2^\cont$.
\end{proof}

\begin{example}\label{ex:extref}
There is a non-reflexive pseudocompact group $G$ that arises
as an extension of a compact group $G/L\cong \Z(2)^{2^\cont}$ by
a closed, reflexive, $\omega$-bounded subgroup $L$ of $G$.
\end{example}

\begin{proof}
Let $P$ be the $P$-modification of the compact group $\Z(2)^{\cont}$.
By Theorem~4.8 of \cite{jlm}, $P$ is  reflexive. Denote by $L$ the
character group of $P$. It is easy to see that $|L|=2^\cont$. By
Lemma~\ref{pdualsigma}, there is a topological monomorphism
$j\colon L\to \Z(2)^{2^\cont}$ such that $j(L)$ is dense in
$\Z(2)^{2^\cont}$ and $j(L)\cap \Sigma\Z(2)^{2^\cont}=\{e\}$.
In the sequel we identify $L$ with $j(L)$.

Let $K$ be a compact subset of $\Z(2)^{2^\cont}$ such that
$|K|=2^\cont$ and $\hull{K}$ is dense in $\Z(2)^{2^\cont}$. For
example, one can take $K=\{e\}\cup\{b_\alpha: \alpha<2^\cont\}$,
where $e$ is the neutral element of $\Z(2)^{2^\cont}$ and
$b_\alpha(\beta)=-1$ only if $\alpha=\beta$; $\alpha,\beta<2^\cont$.
Let also $R_K$ be a subgroup of $\Z(2)^{2^\cont}$ such that
$\Z(2)^{2^\cont}=\la K\ra \oplus R_K$ and define a (necessarily
discontinuous) homomorphism $\phi \colon \Z(2)^{2^\cont}\to \Z(2)^{2^\cont}$
with the following properties:
\begin{enumerate}
  \item $\phi\rest{\la K\ra}$ is the identity mapping;
  \item $\phi(r)\notin rL$ for every $r\in R_K$ with $r\neq e$.
\end{enumerate}
To construct such a homomorphism we first observe that
$|\Z(2)^{2^\cont}/\hull{K}|=2^{2^{\cont}}$ and then consider a maximal
independent subset $\{r_\alpha \colon \alpha<2^{2^\cont}\}$ of $R_K$.
It suffices to define $\phi$ satisfying (1) such that
$\phi(r_\alpha)\notin \la \bigcup\{r_\beta L\cup\phi(r_\gamma)L
\colon \beta\leq \alpha,\: \gamma<\alpha\} \ra$, for each $\alpha<2^\cont$.

We use the homomorphism $\phi$ to apply Lemma~\ref{lem:constG}
and consider the subgroup $G=\Gr(\phi)\cdot N$ of $\Pi=
\Z(2)^{2^\cont}\times\Z(2)^{2^\cont}$, where $N=\{e\}\times L$.
By Lemma~\ref{lem:constG}, $G$ is a $G_\delta$-dense subgroup
of the compact group $\Pi$, so $G$ is pseudocompact according
to \cite[Theorem~1.2]{CR2}.

Then $\widetilde{K}=\{(k,k)\colon k\in K\}$ is a compact subset of
$G$. If $(a,b)\in G\cap cl_{\Pi}\la \widetilde{K}\ra$, then $a=b$ and
$a=\phi(a)y$ with $y \in L$. Since $a=kr$ with $k\in \la K\ra$ and
$r\in R_K$, we see $ry^{-1}=\phi(r)$. Therefore (2) implies
that $r=e$. This proves the inclusion
\[
cl_{G} \la \widetilde{K}\ra\subset
\{(a,a)\colon a\in \la K\ra\},
\]
while the inverse inclusion is evident. Thus $\hull{\widetilde{K}}$
is a closed subgroup of $G$. Since $\hull{K}$ is a proper dense
subgroup of $\Z(2)^{2^\cont}$ it follows that $\la \widetilde{K}\ra$
is not compact and, hence, the group $G$ is not reflexive by
Corollary~\ref{cor:boundedpseudo}.
\end{proof}

We now construct a larger family of extensions of compact (even metrizable)
groups by $\omega$-bounded groups. This requires several preliminary steps.

Given an abstract Abelian group $G$, we denote by $G^\#$ the underlying
group $G$ which carries the maximal precompact group topology \cite{vDo}.
This topology on $G$ is called the \textit{Bohr} topology of $G$. Notice
that \textit{every} homomorphism of $G^\#$ to $\T$ is continuous and that
this property characterizes $G^\#$ among precompact groups. The following
fact is a part of the duality folklore.

\begin{lemma}\label{Le:vDo}
Let $K$ be a compact Abelian topological group and $(K^\wedge,\tau_p(K))$
be the dual group of $K$ with topology $\tau_p(K)$ of pointwise convergence
on elements of $K$. Then $\tau_p(K)$ is the Bohr topology of the (abstract)
group $K^\wedge$.
\end{lemma}

\begin{proof}
By the classical Pontryagin--van~Kampen duality theorem
$K^\wedge$ is a discrete group and $\alpha_K\colon K\to K^{\ww}$
is a topological isomorphism.  Therefore, the Bohr topology of
$K^\wedge$ is the precompact group topology $\tau_p(K)$
generated by $K$.
\end{proof}

\begin{definition}\label{Df1}
Let $\kappa\geq \omega$ be a cardinal number. We say that a subgroup
$L\leq \Z(2)^{2^\kappa}$ satisfies condition $(Sm)$ provided that for
every $N\in[L]^{\leq \omega}$ there is a set $A_N\subset 2^\kappa$
with $|A_N|=2^\kappa$ such that $\pi_{A_N}(N)=\{e_N\}$,
where $\pi_{A_N}\colon \Z(2)^{2^\kappa}\to \Z(2)^{A_N}$ is the
projection and $e_N$ is the neutral element of $\Z(2)^{A_N}$.
\end{definition}

Given a cardinal $\kappa$, we denote by $\Z(2)^{(\kappa)}$ the
direct sum of $\kappa$ copies of the group $\Z(2)$.

\begin{lemma}\label{lem:newref}
Let $\kappa\geq \omega$ be a cardinal, $L$ be a dense pseudocompact
subgroup of $ \Z(2)^{2^\kappa}$ with $|L|\leq 2^\kappa$, and suppose
that $L^\wedge$ is discrete and satisfies $(Sm)$. Then there exists a
pseudocompact group $G$ containing $L$ as a closed subgroup with
$G/L$ compact and such that
\[
G^\wedge=\left(\Z(2)^{(\kappa)} \right)^\# \times L^\wedge,
\]
where $\left(\Z(2)^{(\kappa)}\right)^\#$ stands for the group
$\Z(2)^{(\kappa)}$ equipped with its Bohr topology.
\end{lemma}

\begin{proof}
Let $\mathcal{S}$ be the family of countable subgroups of $\Z(2)^\kappa$.
Consider the family
$$
\mathcal{A}=\left\{(S,f,h)\colon
S\in\mathcal{S},\ f\in Hom(S, \Z(2)),\ h\in Hom(S, L) \right\}.
$$
It is clear that $|\mathcal{A}|=2^\kappa$. Let us fix an injective mapping
$\rho \colon \mathcal{A}\to 2^{\kappa}$ such that $\rho(S,f,h)\in A_{h(S)}$
for every $(S,f,h)\in \mathcal{A}$; here $A_{h(S)}$ is as in
Definition~\ref{Df1}. The existence of such a mapping $\rho$ follows
from the equalities $|\mathcal{A}|=2^\kappa$ and $|A_S|=2^\kappa$
for $S\in\mathcal{S}$. For each $(S,f,h)\in \mathcal{A}$, we consider
a homomorphism $\varphi_\alpha\colon \Z(2)^\kappa\to \Z(2)$ extending
$f$, where $\alpha=\rho(S,f,h)$. Finally, for each  $\alpha<2^\kappa$ we
define a homomorphism $\phi_\alpha\colon \Z(2)^\kappa\to \Z(2)$ by
the following rule:
\[
\phi_\alpha= \begin{cases}   \varphi_\alpha, \mbox{ if }
\alpha=\rho(S,f,h) \mbox{ for some } (S,f,h)\in \mathcal{A};\\
\mathbf{1}  \mbox{ (constant),  if } \alpha \notin
\rho(\mathcal{A}).\end{cases}
\]
Let $\phi$ be the diagonal product of $\phi_\alpha$'s,
${\displaystyle \phi\colon \Z(2)^\kappa\to \Z(2)^{2^\kappa} }$.

We define the group $G=\Gr(\phi)\cdot N$ for these $\phi$ and
$N$ as in Lemma~\ref{lem:constG}. Then $G$ is a dense subgroup
of $\Pi=\Z(2)^\kappa\times\Z(2)^{2^\kappa}$. According to
Lemma~\ref{lem:constG}, $G$ is $G_\delta$-dense in $\Pi$,
so $G$ is pseudocompact by \cite[Theorem~1.2]{CR2}. Let $\pi$
be the projection of the product $\Pi$ to the first factor $\Z(2)^\kappa$.

We now proceed to identify the dual group of $G$. This task requires
an analysis of the structure of compact subsets of $G$. We start
with the following fact:\smallskip

{\bf Claim~1.} \textit{If $S$ is a countable subgroup of $G$ and the
restriction of $\pi$ to $S$ is one-to-one, then $S$ is $h$-embedded
in $G$.}\smallskip

Indeed, let $S$ be a countable subgroup of $G$ such that $\pi\restri_S$
is one-to-one. Then $S$ is the graph of a homomorphism
$g\colon D\to\Z(2)^{2^\kappa}$, where $D=\pi(S)$. It follows from
the definition of $G$ that $g=h\cdot\phi\rest{D}$, where $h$ is a
homomorphism of $D$ to $L$.

Let $p\colon S\to\Z(2)$ be an arbitrary homomorphism. For every
$z\in S$, put $f(\pi(z))=p(z)$. Then $f\colon D\to\Z(2)$ is a homomorphism
and $(D,f,h)\in \mathcal{A}$. With $\alpha=\rho(D,f,h)$ we have
$\phi_\alpha\restr_D\,=\varphi_\alpha\restr_D\,=f$ and $h(x)_\alpha=1$
for each $x\in D$ since $\alpha\in A_{h(D)}$.

For every $\alpha<2^\kappa$, let $\pi_\alpha\colon\Z(2)^{2^\kappa}\to\Z(2)$
be the projection of the product group $\Z(2)^{2^\kappa}$ to the
$\alpha$th factor $\Z(2)_{(\alpha)}$. We claim that the homomorphism
$p$ coincides with the restriction to $S$ of the continuous homomorphism
$\pi_\alpha\circ\varpi$, where $\varpi\colon \Pi\to\Z(2)^{2^\kappa}$
is the projection of $\Pi$ to the second factor. Indeed, take an arbitrary
element $z=(x,y)\in S$. Then $y=h(x)\cdot\phi(x)$, since $S$ is the graph
of the homomorphism $h\cdot\phi\rest_D$. We have, on one side, that
\begin{equation}\label{1}
p(z)=f(\pi(z))=f(x)=\varphi_{\alpha}(x).
\end{equation}
On the other side, it follows from our choice of $\alpha<2^\kappa$ and
the definition of $\varphi$ that
\begin{eqnarray*}
\pi_\alpha(\varpi(z))=\pi_\alpha(y)\hskip-4pt &=&\hskip-4pt
\pi_\alpha(h(x))\cdot\pi_\alpha(\phi(x))\\
\hskip-4pt &=&\hskip-4pt h(x)_\alpha\cdot\varphi_{\alpha}(x)
=\varphi_{\alpha}(x).
\end{eqnarray*}
Comparing the above equality and equality~\eqref{1}, we infer that
$p(z)= \pi_\alpha(\varpi(z))$ for each $z\in S$, which proves that
$p=\pi_\alpha\circ\varpi\rest{S}$.

Thus, every homomorphism $p\colon S\to\Z(2)$ extends to
a continuous homomorphism of $G$ to $\Z(2)$. This proves
Claim~1.\smallskip

{\bf Claim~2.} \textit{If $S$ is a countable subgroup of $G$ and
the restriction of $\pi$ to $S$ is one-to-one, then
$\pi(cl_G(S))=\pi(S)$.}\smallskip

Let $g$ be an element of $G$ such that $\pi(g)\notin\pi(S)$. It suffices
to show that $g\notin cl_G(S)$. It follows from our choice of $g$ that
the restriction of $\pi$ to the subgroup $T=S\cdot \hull{g}$ of $G$ is also
one-to-one. Hence, by Claim~1, $T$ is $h$-embedded in $G$. Let
$\varphi$ be a homomorphism of $T$ to $\Z(2)$ such that
$\varphi(g)=-1$ and $\varphi(S)=\{1\}$. Then $\varphi$ is continuous
on $T$ and, therefore, $g\notin cl_G(S)$.\smallskip

{\bf Claim~3.} \textit{The projection $\pi(K)$ is finite, for every compact
set $K\subset G$.}\smallskip

Suppose for a contradiction that $\pi(K)$ is infinite, for a compact
subset $K$ of $G$. We can assume without loss of generality that
$\pi(K)$ does not contain the neutral element $e$ of $\Z(2)^\kappa$.
Clearly $\pi(K)$ contains a countable infinite independent subset,
say $X$. Choose a subset $Y$ of $K$ such that $\pi(Y)=X$ and
the restriction of $\pi$ to $Y$ is one-to-one. Then $Y$ is countable
and independent in $G$ and the restriction of $\pi$ to the subgroup
$S=\hull{Y}$ of $G$ is one-to-one. Let $C=K\cap cl_G(S)$.
Then $C$ is a compact subset of $G$ and Claim~2 implies that
$\pi(C)\subset\pi(S)$. It also follows from $Y\subset K\cap S\subset C$
and $\pi(Y)=X$ that the compact sets $C$ and $\pi(C)$ are infinite.
Since $\pi(C)$ is countable (hence metrizable) and $X\subset\pi(C)$,
there exists a sequence $\{x_n: n\in\omega\}\subset X$ converging
to an element $x^*\in\pi(C)$, where $x^*\neq x_n$ for each $n\in\omega$.
Notice that $x^*\neq e$. By induction we can choose an infinite subset
$X'$ of $\{x_n: n\in\omega\}$ such that $x^*\notin\hull{X'}$. Take a
subset $Y'$ of $Y$ such that $\pi(Y')=X'$ and let $S'=\hull{Y'}$.
Arguing as above, we see that $C'=K\cap cl_G(S')$ is a compact
subset of $G$, $\pi(C')\subset \pi(S')$, and that $X'=\pi(Y')\subset \pi(C')$.
In particular, the compact set $\pi(C')$ contains infinitely many points
$x_n$'s and, hence, $x^*\in\pi(C')$. The latter, however, is impossible
since $\pi(C')\subset\pi(S')=\hull{X'}$ and $x^*\notin\hull{X'}$. This
proves Claim~3.\smallskip

We now obtain a complete description of $G^\wedge$, both
algebraic and topological. Since $G$ is dense in $\Pi$ (and
hence each character of $G$ extends to a character of
$\Pi$), there exists a natural (abstract) isomorphism
\begin{equation}
G^\wedge\cong(\Z(2)^\kappa\times\Z(2)^{2^\kappa})^\wedge
\cong(\Z(2)^\kappa)^\wedge
\oplus (\Z(2)^{2^\kappa})^\wedge.\label{Eq:Iso}
\end{equation}
The second isomorphism in (\ref{Eq:Iso}) is obtained by restricting
every character of $\Pi$ to the factors $\Z(2)^\kappa$ and $\Z(2)^{2^\kappa}$.
Since the groups $\Z(2)^\kappa$ and $\Z(2)^{2^\kappa}$ are compact,
$(\Z(2)^\kappa)^\wedge$ and $(\Z(2)^{2^\kappa})^\wedge$ are Boolean
groups of cardinality $\kappa$ and $2^\kappa$, respectively.

Finally, we claim that $G^\wedge$ is topologically isomorphic to the group
\begin{equation}
((\Z(2)^\kappa)^\wedge,\tau_p) \oplus L^\wedge,\label{Eq:TopIso}
\end{equation}
where $\tau_p$ stands for the pointwise convergence topology on
the abstract group $(\Z(2)^\kappa)^\wedge$.

For each $F\subset \Z(2)^\kappa$ and each $P\subset G$, let
$$
C_{F,P}=\{(x,\phi(x)y)\colon x\in F,\; y\in P \}.
$$
If $K\subset G$ is compact, then $F=\pi(K)$ is finite and the set
$$
P_K=\{y\in L\colon (x,\phi(x)y)\in K \mbox{ for some } x\in F\}
$$
is a compact subset of $L$. Since $K\subset C_{F,P_K}$, we
deduce that the family
\begin{equation}
\{(C_{F,P})^\vart: F\subset \Z(2)^\kappa,\ |F|<\omega,\
P\subset L,\ P \mbox{ is compact}\}\label{Eq:8}
\end{equation}
forms a local base at the neutral element of the group $G^\wedge$.

Let now $F\subset \pi(G)$ be finite and $P \subset L$ be compact.
If $\widetilde{P}=\phi(F)P$, then
\[
F^\vart \times \widetilde{P}^\vart \subset (C_{F,P})^\vart.
\]
The sets $(C_{F,P})^\vart$ are therefore neighborhoods of the
identity in $((\Z(2)^\kappa)^\wedge,\tau_p) \times L^\wedge$.

We now use the fact that $L^\wedge$ is discrete and take a compact
set $P_0\subset L$ such that $P_0^\vart =\{\mathbf{1}\}$, where
$\mathbf{1}$ is the neutral element of $L^\wedge$. Since $C_{F,P_0}$
contains the set $\{e\}\times P_0$, we see that
$F^\vart\times\{e\}=(C_{F,P_0})^\vart$, for every finite set
$F\subset\Z(2)^\kappa$. Therefore, $G^\wedge$ is topologically
isomorphic to $((\Z(2)^\kappa)^\wedge,\tau_p) \times L^\wedge$.
It only remains to observe that, by Lemma~\ref{Le:vDo},
$((\Z(2)^\kappa)^\wedge,\tau_p)$ is exactly $\left(\Z(2)^{(\kappa )}\right)^\#$.
\end{proof}

\begin{corollary}\label{Co:OldT}
There exists a pseudocompact Abelian group $G$ which contains a
closed $\omega$-bounded (hence countably compact) subgroup $N$
such that the dual groups $N^\wedge$ and $(G/N)^\wedge$ are discrete,
but $G^\wedge$ is not. In addition, the quotient group $G/N$ is compact
metrizable, while the bidual group $G^{\wedge\wedge}$ is compact and
topologically isomorphic to $\varrho{G}$, the completion of $G$.
\end{corollary}

\begin{proof}
In Lemma~\ref{lem:newref}, let $\kappa=\omega$ and take $L$
to be the $\Sigma$-product $\Sigma\Pi$, where  $\Pi=\Z(2)^{2^\cont}$.
It is obvious that $L$ is $\omega$-bounded and dense in $\Pi$, satisfies
condition $(Sm)$, and the dual group $L^\wedge$ is discrete.
Applying Lemma~\ref{lem:newref}, we find a dense pseudocompact
subgroup $G$ of $\Z(2)^\omega\times\Pi$ containing $L$ as a closed
subgroup such that $G/L\cong \Z(2)^\omega$ and
$G^\wedge\cong (\Z(2)^{(\omega)})^\#\times L^\wedge$. In particular,
$G^\wedge$ is not discrete. Further, the standard calculation shows that
 \[
 G^{\wedge\,\wedge}=\Z(2)^{\omega}\times L^{\wedge\wedge}=
 \Z(2)^\omega \times \Z(2)^{2^\cont}.
 \]
Since $G$ is dense in $\Z(2)^\omega \times \Z(2)^{2^\cont}$, we
conclude that $G^{\wedge\,\wedge}\cong\varrho{G}$.
\end{proof}

\begin{remark}
\emph{The group $G$ in Corollary~\ref{Co:OldT} fails to be countably compact,
even when the closed subgroup $L$ of $G$ is countably compact (even
$\omega$-bounded) and the quotient group $G/L\cong\Z(2)^\omega$ is
compact and metrizable. The first example of a such a group was
constructed in \cite{BT}.}

Let us show that the present group $G$ contains an infinite closed
discrete subset. We keep the notation adopted in Lemma~\ref{lem:newref}.
Take a sequence $\{x_n: n\in\omega\}\subset \Z(2)^\omega$
converging to an element $x^*\in\Z(2)^\omega$, where $x^*\neq e$.
We can assume that for each $n\in\omega$, the element $x_n$
is not in the subgroup of $\Z(2)^\omega$ generated by the set
$\{x^*\}\cup \{x_k: k<n\}$. For every $n\in\omega$, let $z_n=(x_n,y_n)$,
where $y_n=\phi(x_n)$. It follows from the definition of $G$ that
the set $P=\{z_n: n\in\omega\}$ is contained in $G$.

We claim that $P$ is closed and discrete in $G$. Since
$\pi(z_n)=x_n$ and $x_n\to x^*$, all accumulation points of $P$,
if any, lie in $\pi^{-1}(x^*)$. Take an arbitrary point $z=(x^*,y)\in
\pi^{-1}(x^*)\cap G$. Again, our definition of $G$ implies that
$y=\phi(x^*)\cdot{s}$, for some $s\in L$. Denote by $D$ the
subgroup of $\Z(2)^\omega$ generated by $\{x^*\}\cup\{x_n: n\in\omega\}$
and take a homomorphism $f\colon D\to\Z(2)$ such that $f(x^*)=-1$
and $f(x_n)=1$, for each $n\in\omega$. Since $\supp(s)\subset\cont$ is
countable and the mapping $\varrho\colon \mathcal{A}\to\cont$ is injective,
there exists a homomorphism $h\colon D\to L$ such that $\alpha=\varrho(D,f,h)
\notin\supp(s)$. Notice that $\varphi_\alpha\restri_D=f$. We now have that
\begin{eqnarray*}
\pi_\alpha\varpi(z)=\pi_\alpha(y) &=&
\pi_\alpha(\phi(x^*)\cdot{s})\\
& = &
\pi_\alpha\phi(x^*)\cdot\pi_\alpha(s)=\varphi_\alpha(x^*)=f(x^*)=-1,
\end{eqnarray*}
while a similar calculation shows that $\pi_\alpha\varpi(z_n)=1$,
for each $n\in\omega$. Since the homomorphism $\pi_\alpha\circ\varpi
\colon \Z(2)^\omega\times\Z(2)^\cont\to\Z(2)_{(\alpha)}$ is continuous,
we conclude that $z\notin cl_G{P}$. This proves that $P$ is closed and
discrete in $G$ and that $G$ is not countably compact.
\end{remark}

\bibliographystyle{plain}

\begin{thebibliography}{99}

\def\cprime{$'$} \def\cprime{$'$} \def\cprime{$'$} \def\cprime{$'$}
  \def\polhk#1{\setbox0=\hbox{#1}{\ooalign{\hidewidth
  \lower1.5ex\hbox{`}\hidewidth\crcr\unhbox0}}}
  \def\polhk#1{\setbox0=\hbox{#1}{\ooalign{\hidewidth
  \lower1.5ex\hbox{`}\hidewidth\crcr\unhbox0}}} \def\cprime{$'$}

\bibitem{ardaetal}
S.~Ardanza-Trevijano, M.\,J.~Chasco, X.~Dom\'{\i}nguez, and
M.\,G.~Tkachenko,
\newblock Precompact noncompact reflexive abelian groups,
\newblock \textit{Forum Math.} {\bf 24} (2) (2012), 289--302.

\bibitem{AT} A.\,V.~Arhangel'skii and M.\,G.~Tkachenko,
\newblock \textit{Topological Groups and Related Structures},
\newblock Atlantis Series in Mathematics, Vol.~I,
Atlantis Press and World Scientific, Paris--Amsterdam 2008.

\bibitem{Aub} L.~Au{\ss}enhofer,
\newblock \emph{Contributions to the Duality theory of Abelian
topological groups and to the theory of nuclear groups},
\newblock Ph.D. Dissertation (Der Mathematischen Fakult\"at der
Eberhard-Karls-Universit\"at T\"ubingen zur Erlangung des Grades
eines Doktors der Naturwissenscheaften), 1998.

\bibitem{bana91}
W.~Banaszczyk, \emph{Additive subgroups of topological vector
spaces}, Springer-Verlag, Berlin, 1991.

\bibitem{BCM} W.~Banaszczyk, M.\,J.~Chasco, and
E.~Mart\'{\i}n-Peinador,
\newblock Open subgroups and Pontryagin duality,
\newblock  \textit{Math. Z.} \textbf{215} (2) (1994), 195--204.

\bibitem{BT} M.~Bruguera and M.~Tkachenko,
\newblock Extensions of topological groups do not respect countable
compactness,
\newblock \textit{Questions and Answers in Gen. Topol.}
\textbf{22} (1) (2004), 33--37.

\bibitem{BT2} M.~Bruguera and M.~Tkachenko,
\newblock The three space problem in topological groups,
\newblock \textit{Topology Appl.} \textbf{153} (3) (2006), 2278--2302.

\bibitem{Cha} M.\,J.~Chasco,
\newblock Pontryagin duality for metrizable groups,
\newblock \textit{Arch. Math.} \textbf{70} (1) (1998), 22--28.

\bibitem{comfgali03} W.\,W.~Comfort and J.~Galindo,
\newblock Pseudocompact topological group refinements of maximal weight,
\newblock \textit{Proc. Amer. Math. Soc.} \textbf{131} (4) (2003), 1311--1320.


\bibitem{comfraczktrig04} W.\,W.~Comfort and S.\,U.~Raczkowski and F. \,J.~Trigos-Arrieta
\newblock The dual group of a dense subgroup
\newblock \textit{Czechoslovak Math. J.} \textbf{54} (2004), 509--533.

\bibitem{CR} W.\,W.~Comfort and K.\,A.~Ross,
\newblock Topologies induced by groups of character,
\newblock \textit{Fund. Math.} \textbf{55} (1964), 283--291.

\bibitem{CR2} W.\,W.~Comfort and K.\,A.~Ross,
\newblock Pseudocompactness and uniform continuity in topological groups,
\newblock \textit{Pacific J. Math.} \textbf{16} (1966), 483--496.

\bibitem{DS} S.~Dierolf and U.~Schwanengel,
\newblock Examples of locally compact non-compact minimal
topological groups,
\newblock \textit{Pacific J. Math.} \textbf{82} (1979), 349--355.

\bibitem{vDo} E.\,K.~Douwen,~van,
\newblock The maximal totally bounded group topology on $G$
and the biggest minimal $G$-space for Abelian groups $G$,
\newblock \textit{Topol. Appl.} \textbf{34} (1990), 69--91.

\bibitem{GH} J.~Galindo and S.~Hern\'andez,
\newblock Pontryagin--van Kampen reflexivity for free Abelian topological
groups,
\newblock \textit{Forum Math.} \textbf{11} (4) (1999), 399--415.

\bibitem{galimaca11} J.~Galindo and S.~Macario,
\newblock Pseudocompact group topologies with no infinite compact subsets,
\newblock {\em J. Pure Appl. Algebra} \textbf{215} (4) (2011), 655--663.

\bibitem{jlm}
J.~Galindo, L.~Recoder-N\'u\~{n}ez, M.~Tkachenko,
\newblock Nondiscrete $P$-groups can be reflexive,
\newblock {\em Topology Appl.} \textbf{158} (2) (2011), 194--203.

\bibitem{jlm2}
J.~Galindo, L.~Recoder-N\'u\~{n}ez, M.~Tkachenko,
\newblock Reflexivity of prodiscrete topological groups,
\newblock {\em J. Math. Anal. Appl. } \textbf{384} (2) (2011), 320--330.

\bibitem{gilljeri}
L.~Gillman and M.~Jerison,
\newblock \textit{Rings of Continuous Functions},
\newblock Van Nostrand, New York, 1960.

\bibitem{Gra} D.\,L.~Grant,
\newblock Topological groups which satisfy an open mapping theorem,
\newblock \textit{Pacific J. Math.} \textbf{68} (1977), 411--423.

\bibitem{hernmaca03}
S.~Hern\'andez and S.~Macario,
\newblock Dual properties in totally bounded abelian groups,
\newblock {\em Arch. Math. (Basel)} \textbf{ 80} (3) (2003), 271--283.

\bibitem{hernuspe00}
S.~Hern\'andez and V.~Uspenskij,
\newblock Pontryagin duality for spaces of continuous functions,
\newblock {\em J. Math. Anal. Appl.} \textbf{242} (2) (2000), 135--144.

\bibitem{Nob} N.~Noble,
\newblock Products with closed projections,
\newblock \textit{Trans. Amer. Math. Soc.} \textbf{140} (1969), 381--391.

\bibitem{Nob2} N.~Noble,
\newblock $k$-groups and duality,
\newblock \textit{Trans. Amer. Math. Soc.} \textbf{151} (1970), 551--561.

\bibitem{Pes} V.\,G.~Pestov,
\newblock Free topological abelian groups and the Pontryagin duality,
\newblock \textit{Mosc. Univ. Math. Bull.} \textbf{41} (1) (1986), 1--4.
\newblock Russian original in: \textit{Vestnik Mosk. Univ. Ser. I\/} no.~1
(1986), 3--5.

\bibitem{Pes2} V.\,G.~Pestov,
\newblock Free abelian topological groups and the Pontryagin--van Kampen
duality,
\newblock \textit{Bull. Austral. Math. Soc.} \textbf{52} (1995), 297--311.

\end{thebibliography}
\def\cprime{$'$} \def\cprime{$'$} \def\cprime{$'$} \def\cprime{$'$}
  \def\polhk#1{\setbox0=\hbox{#1}{\ooalign{\hidewidth
  \lower1.5ex\hbox{`}\hidewidth\crcr\unhbox0}}}
  \def\polhk#1{\setbox0=\hbox{#1}{\ooalign{\hidewidth
  \lower1.5ex\hbox{`}\hidewidth\crcr\unhbox0}}} \def\cprime{$'$}

\section{P-modification of reflexive groups}\label{Sec:Pmod}
In our first example we show that the $P$-modification of a reflexive
$\sigma$-compact group can fail to be reflexive. Our argument uses
essentially Pestov's theorem about the reflexivity of free Abelian
topological groups on zero-dimensional compact spaces (see
\cite{Pes, Pes2}).

Let $D$ be an uncountable discrete space and $X$ a one-point
compactification of $D$ with a single non-isolated point $x_0$.

\begin{lemma}\label{Le:1}
Let $G=A(X)$ be the free Abelian topological group over $X$.
Then the $P$-modification $PG$ of $G$ is topologically isomorphic
to the free Abelian topological group $A(Y)$, where $Y=PX$ is the
$P$-modification of the space $X$.
\end{lemma}

\begin{proof}
It is clear that $Y$ is a Lindel\"of $P$-space, and so is every finite
power of $Y$ \cite{Nob}. Hence the group $A(Y)$ is Lindel\"of, while
\cite[Proposition~7.4.7]{AT} implies that $A(Y)$ is a $P$-space.

Let $i\colon A(Y)\to A(X)$ be the continuous isomorphism of $A(Y)$
onto $A(X)$ which extends the identity mapping of $Y$ onto $X$.
It suffices to verify that $i$ is a homeomorphism of $A(Y)$ onto
$PA(X)$ (the group $A(X)$ with the $P$-modified topology).

Let $C$ be a countable subset of $D$. Denote by $r_C$ the retraction
of $X$ onto $X_C=C\cup\{x_0\}$, where $r_C(x)=x$ for each $x\in X_C$
and $r_C(y)=x_0$ for each $y\in X\setminus X_C$. Clearly, $r_C$ is continuous.
Extend $r_C$ to a continuous homomorphism $R_C\colon A(X)\to A(X_C)$.
Since the free Abelian topological group $A(X_C)$ has countable
pseudocharacter (and $R_C$ is continuous), $\ker R_C$ is a closed
$G_\delta$-set in $A(X)$. Hence $H_C=i^{-1}(\ker{R_C})$ is an open
subgroup of $A(Y)$.

Consider the family
\[
\mathcal{H}=\{H_C: C\subset D,\ |C|\leq\omega\}.
\]
It is easy to see that if $\{C_n: n\in\omega\}$ is a sequence of countable
subsets of $D$ and $C=\bigcup_{n\in\omega} C_n$, then $H_C\subset
\bigcap_{n\in\omega} H_{C_n}$. It is also clear that the intersection of
the family $\mathcal{H}$ contains only the neutral element of $A(Y)$.

We claim that $\mathcal{H}$ is a local base at the neutral element $e$
of $A(Y)$. Indeed, take an arbitrary open neighborhood $U$ of $e$ in
$A(Y)$. Suppose for a contradiction that every element of $\mathcal{H}$
meets the closed subset $F=A(Y)\setminus U$ of $A(Y)$. Then the family
$\mathcal{H}\cup\{F\}$ of closed subsets of $A(Y)$ has the countable
intersection property. Since the space $A(Y)$ is Lindel\"of, we must
have that $F\cap\bigcap\mathcal{H}\neq\emptyset$, which contradicts
the equality $\bigcap\mathcal{H}=\{e\}$.

Finally, since $\ker{R_C}$ is open in $PA(X)$,
the $P$-modification of $A(X)$, this implies that $i\colon A(Y)\to PA(X)$
is a homeomorphism.
\end{proof}

The proof of the following fact is implicitly contained in (the proof of)
Lemma~17.11 of \cite{bana91}. We include the proof for the sake of
completeness, it follows that the same result is valid for any infinite
direct sum or any infinite $\Sigma$-product of compact groups (see \cite[Lemma 3.12]{comfraczktrig04})

\begin{lemma}\label{sigma}
The dual group $\left(\Sigma\T^D\right)^\wedge$ of the $\Sigma$-product
$\Sigma\T^D$ is  a discrete group.
\end{lemma}

\begin{proof}
Define for each $x\in \T$ and $d\in D$ the element of $f_{x,d}\in \T^D$
that takes the value $x$ at $d$ and $1$ elsewhere.
The set $K=\{f_{x,d}\colon x\in \T,\ d\in D\}$ is then a compact subset of
$\Sigma\T^D$.  Since for each $d\in D$, the set $A_d=\{f_{x,d} \colon x\in \T\}$
is a subgroup of $\T^D$ contained in $K$, we conclude that $\chi(A_d)=\{1\}$
for every $\chi \in K^\rhd$ and $d\in D$. But  the subgroup generated by all
the $A_d$'s is  dense in $\T^D$, therefore $K^\rhd=\{1\}$ and
$\left(\Sigma\T^D\right)^\wedge$ is discrete.
\end{proof}

\begin{example}\label{Ex:1}
{\it The free Abelian topological group $G=A(X)$ is reflexive, but
the $P$-modification $PG$ of $G$ is not. Furthermore, the second
dual of $PG$ is discrete.}
\end{example}

\begin{proof}
Since $X$ is a zero-dimensional compact space, the reflexivity of $G$
follows from \cite{Pes} or \cite{Pes2}. Let us verify that $PG$ is not reflexive.
Since, by Lemma~\ref{Le:1}, $PA(X)$ is topologically isomorphic to $A(Y)$,
where $Y=PX$, and neither $Y$ nor $A(Y)$ is discrete, it suffices to show
that the second dual of $A(Y)$ is discrete.

Clearly, all compact subsets of the $P$-space $Y$ are finite, so $Y$ is
a $\mu$-space. It follows from \cite[Theorem~2.1]{GH} that the dual
group $A(Y)^\wedge$ is topologically isomorphic to $C_p(X,\T)$,
the group of all continuous functions on $X$ with values in the circle
group $\T$, endowed with the pointwise convergence topology.

Let $Y=D\cup\{x_0\}$, where $x_0$ is the only non-isolated point of
$Y$. Every neighborhood of $x_0$ in $Y$ has the form $Y\setminus C$, where
$C$ is a countable subset of $D$. Therefore, for every element $f\in C_p(Y,\T)$,
there exists a countable set $C\subset D$ such that $f(x)=f(x_0)$ for each
$x\in Y\setminus C$.

Denote by $\Sigma\T^D$ the $\Sigma$-product lying in $\T^D$ and
considered as a dense subgroup of the compact group $\T^D$.
We consider a mapping $\varphi\colon C_p(Y,\T)\to \T\times\Sigma\T^D$
defined by $\varphi(f)=(f(x_0),t_f\cdot{f})$, where $t_f=f(x_0)^{-1}\in\T$
and the function $t_f\cdot{f}$ is restricted to $D$. Then $t_f\cdot{f}\in
\Sigma\T^D$ and $\varphi(f)\in\T\times\Sigma\T^D$. Since $C_p(Y,\T)$
carries the topology of pointwise convergence, $\varphi$ is a topological
isomorphism of $C_p(Y,\T)$ onto $\T\times\Sigma\T^D$. Hence
the dual of $C_p(Y,\T)$ is topologically isomorphic to
$(\T\times\Sigma\T^D)^\wedge\cong \T^\wedge\times (\Sigma\T^D)^\wedge
\cong\Z_d\times (\Sigma\T^D)^\wedge$, where $\Z_d$ is the discrete group
of integers. Finally, we know by Lemma~\ref{sigma} that the dual group
$(\Sigma\T^D)^\wedge$ is discrete. Hence the second dual
$A(Y)^{\wedge\wedge}$ is discrete as well.
\end{proof}

\begin{corollary}\label{Cor:Last}
Let $D$ be an uncountable discrete space and let $Y$ denote the one-point
Lindel\"ofication of $D$. Then $C_p(Y,\T)$ is not reflexive.
\end{corollary}

\begin{proof}
We note that $Y=PX$, where $X$ is the one-point compactification
of $D$. The proof of Example \ref{Ex:1} shows that $C_p(Y,\T)^\wedge$ is
discrete, while $C_p(Y,\T)$, being  a proper dense subgroup of $\T^Y$, is
not compact.
\end{proof}

\begin{remark}
Corollary~\ref{Cor:Last} is in  contrast with Example~3.12 of
\cite{hernuspe00}, where it is shown that $C_p(Y,\C)$ is reflexive
if  $|D|=\omega_1$ and MA($\omega_1$) is assumed (here $\C$
stands for the field of complex numbers).
\end{remark}

\end{document}